\documentclass[11pt]{article}
 
\usepackage{amsmath, psfrag, epsfig, amsfonts, amsthm}
\newtheorem{theorem}{Theorem}

\newtheorem{remark}{Remark}
\newtheorem{remarks}[remark]{Remarks}
\def\states{\varOmega}
\def\Nrv{\varPhi}
\def\Nrvhat{\widehat\varPhi}
\def\reals{\mathbb{R}}
\def\nats{\mathbb{N}}

\def\Proh{\varrho}
\def\slack{\lambda}
\def\Prohp{\varrho_\slack}

\def\statesP{\widehat\varOmega}
\def\Yhat{\widehat Y}
\def\Uhat{\widehat U}
\def\Shat{\widehat S}
\def\What{\widehat W}
\def\Kbar{\overline K}
\def\PP{\widehat P}
\def\eps{\varepsilon}
\def\defeq{:=}
\def\law{\mathcal{L}}
\def\Xhat{\widehat X}
\def\xbar{\bar x}
\def\ybar{\bar y}
\def\bhat{\hat b}
\def\calE{\mathcal{E}}
\def\unit{u_0}

\begin{document}

\title{\bf On the approximation of one Markov chain by another}

\author{Mark Jerrum\thanks{Partially supported
by EPSRC grant ``Sharper Analysis of Randomized Algorithms:
a Computational Approach'' and the IST Programme of the EU under
contract IST-1999-14036 (RAND-APX)\null.  The work described here 
was partially carried out while the author was visiting 
the Isaac Newton Institute for Mathematical Sciences, Cambridge, UK\null.
Postal address:  School of Informatics, University of Edinburgh, 
The King's Buildings, Edinburgh EH9~3JZ, United Kingdom.}\\
School of Informatics\\
University of Edinburgh
}

\maketitle
 
\begin{abstract}\noindent
Motivated by applications in Markov chain Monte Carlo,
we discuss what it means for one Markov chain to 
be an approximation to another.  Specifically included 
in that discussion 
are situations in which a Markov chain with continuous 
state space is approximated by one with finite state space.
A simple sufficient condition for close approximation 
is derived, which indicates the existence of
three distinct approximation regimes.  
Counterexamples are presented to 
show that these regimes are real and not artifacts of 
the proof technique.  An application to the ``ball walk''
of Lov\'asz and Simonovits is provided as an illustrative
example.
\end{abstract}
 
\newpage
\section{Discussion}
Monte Carlo algorithms compute approximate solutions
to hard problems by extracting information from 
random samples.  Markov chain Monte Carlo (MCMC) algorithms
add an additional ingredient, namely Markov chain simulation,
to this recipe.  The idea is to
devise a Markov chain $(X_t:t\in\nats)$ whose stationary distribution
is the one from which we would like to sample.
The required samples are drawn from 
a realisation of this Markov chain obtained 
by computer simulation. To avoid excessive bias, the 
samples must come from a time step of the realisation
that is beyond the {\it mixing time\/} of the 
Markov chain, i.e., the time~$\tau$ at which $X_\tau$
is close enough to stationarity.

The analysis of MCMC algorithms clearly 
requires us to bound the mixing time from above, 
and several approaches have been proposed 
for achieving this goal.  However, the computer simulation
of the Markov chain will in general be imperfect.
The transition probabilities may not be exactly what 
they should be.  Even worse, the state space may be 
uncountably infinite, so we cannot even represent the
states exactly in the computer.  Does this matter?
Obviously, the answer depends on the accuracy with
which the Markov chain is simulated.  The aim of this 
note is to quantify the required accuracy. 

As a paradigmatic example, consider
the (lazy) ball walk in a convex body, 
due to Lov\'asz and Simonovits~\cite{LS93}.
The state space in this instance
is a convex body in~$\reals^n$, i.e., 
a compact convex set $K\subset\reals^n$ of full dimension.
The transition kernel of the ball walk (with 
step-size $r>0$) is defined 
by the following trial:  Suppose $X_t=x$.  Choose a
point~$y$ uniformly at random (u.a.r.\null) from the ball 
of radius~$r$ centred at~$x$.  If $y\in K$ then
$X_{t+1}$ is~$y$, otherwise $X_{t+1}$ is~$x$. 
The state space is continuous, and the transition kernel also. 
In any implementation it would be necessary
to approximate the states in the realisation of the
ball walk by vectors of finite-precision real numbers;
likewise, the transition kernel would need to 
approximated by some discrete distribution.

This example motivates our general setting.
There is an ``ideal'' ergodic Markov
chain $(\states,P)$, with state space~$\states$ and 
transition kernel~$P$, whose stationary 
distribution and mixing time is known.  Then there 
is a perturbed Markov chain 
$(\statesP, \PP)$, which is the one actually 
implemented.  We assume $\statesP\subseteq\states$. 
Usually, $\statesP$ will be finite, though we don't assume this.
Sometimes, as we have seen, $\states$ will be uncountably 
infinite.
We no not assume that $(\statesP, \PP)$ is necessarily ergodic.
For example, in an implementation of the ball walk,
the low order bits in the finite real number approximations 
might depend deterministically on those of the start state.
In order to compare the $t$-step distributions of the two Markov chains,
we regard $\PP^t(x,\cdot)$ as a probability distribution on~$\states$,
using the convention $\PP^t(x,A)\defeq\PP^t(x,A\cap\statesP)$.

Observe that in general $\PP^t(x,\cdot)$ does not converge 
to $P^t(x,\cdot)$ in usual total variation distance 
(half $\ell_1$-norm).  Indeed, any finite approximation
to the the ball walk will necessarily remain at 
total variation distance~1 throughout,
since $\statesP$ has measure zero in~$\states$.
It is clear, then, that any discussion of finite approximations
to the ball walk must necessarily involve some
underlying metric~$d$ on~$\states$.  In the case of the 
ball walk it would be natural to take $d$ to be the
Euclidean metric.
 
So regard $(\states,d)$ as a metric space, and look at 
convergence in {\it Prohorov metric}:
for Borel probability measures $\pi$ and~$\pi'$ on~$\states$, 
define 
\begin{equation}\label{eq:Prohorov}
\Proh(\pi,\pi')\defeq\inf\big\{\eps:\pi(A)\leq \pi'(A^\eps)+\eps
   \text{ for all closed $A$}\big\},
\end{equation}
where 
$A^\eps\defeq\{y:d(y,A)\leq\eps\}$ and $d(y,A)\defeq\inf\{d(y,x):x\in A\}$.
(It can be shown that taking an infimum just over closed sets
is equivalent to taking an infimum over all Borel sets.)
The appearance of the Prohorov metric in this context is not novel,
as it has been used by a few people, 
for example Diamond et al.~\cite{DPK94},
in studying approximations to dynamical systems.
For reasons that will be mentioned in passing at the relevant moment,
we need the technical condition that $(\states,d)$ is a 
{\it separable\/} metric space.  
This will always be the case in practice
(e.g., for Euclidean space $(\reals^n,\ell_2)$).

Upon reflection, 
there seem to be three prerequisites for $(\statesP, \PP)$ to 
behave as a close approximation to $(\states, P)$.

\begin{enumerate}
\item $\PP(x,\cdot)$ should be close to $P(x,\cdot)$ 
for all $x\in\statesP$.  This is the most obviously necessary 
condition.  The ball-walk example suggests that ``close''
should be measured in the Prohorov metric, and not total variation.

\item $P(x,\cdot)$ should vary smoothly with $x$.  This 
condition is necessary to exclude ``chaotic'' systems 
whose stationary distribution
is very sensitive to small changes in~$P$.

\item $(\states,P)$ should be rapidly mixing.  
Otherwise $(\states,P)$ and $(\statesP,\PP)$ might diverge 
slowly over time, even if conditions (1) and~(2) are met. 
Consider, e.g., a random walk on $\{0,1,\ldots,2^n-1\}$
with a drift of order $2^{-n}$.
\end{enumerate}

Conditions (1) and (3) were noted by Azar et al.~\cite{ABKLP96},
whose motivation was similar to ours, but who considered
the more restricted situation $\statesP=\states$.
They had no need of~(2) since they were dealing only with 
Markov chains with discrete state spaces.

Aside from Azar et al., there is also related work on 
the simulation of dynamical systems, for example,
by Shardlow and Stuart~\cite{ShSt00}.  
Here, the dynamical system may be in continuous time, 
and any computer simulation will involve discretisation
of time as well as of the state space.  (Indeed, it is 
fair to say that the discretisation of time is a greater
concern in this setting.)
Where this work diverges from that in the dynamical systems 
literature is in the emphasis on non-asymptotic bounds 
that explore the dependence of errors on some measure
of the {\it size\/} or {\it complexity\/} 
of the Markov chain.  For example, in the 
simple random walk example from condition~(3) above, 
we are interested in quantifying, in terms of
the size of the state space of the random walk, how close 
the transition kernel $\PP(x,\cdot)$ must be to $P(x,\cdot)$ 
to achieve an adequate approximation.
In the case of the ball walk, we may want to quantify the 
closeness of approximation in terms of the dimension~$n$,
step-size~$r$, and the diameter of the convex body~$K$.
This concern seems less of an 
issue in the dynamical systems literature.

Although Theorem~\ref{thm:main} is billed as the main result,
it must be admitted that its conclusion is unsurprising 
and its proof banal.  Nevertheless, it may 
have some utility in justifying the use of theoretical 
mixing-time upper bounds in imperfect computer simulations,
where real numbers are carried to bounded accuracy
and random variables are sampled from not quite the right 
distributions.  An example application is given in~\S\ref{sec:app}.
The main theoretical contribution of this note is in \S\ref{sec:ex}
where it is shown, through a sequence of counterexamples, that
the three possible behaviours described in Theorem~\ref{thm:main} 
are real, and not artifacts of the proof.  These examples
will hopefully shed light on the main mechanisms at work in 
this setting.

\section{Definitions and preliminaries}
Observe that the two occurrences of $\eps$ in definition~(\ref{eq:Prohorov})
have different functions:  one limits variation in position, 
and the other variation in probability.
In questions of asymptotic convergence it is fine to lump these 
together.  In quantitative work, we want to separate them,
since we need to establish greater control over the former
than the latter.  
In light of this, define a parametric
version of the Prohorov metric
\begin{equation}\label{eq:ParamProhorov}
\Prohp(\pi,\pi')\defeq\inf\big\{\eps:\pi(A)\leq \pi'(A^{\slack\eps})+\eps
   \text{ for all closed $A$}\big\}.
\end{equation}
A metric such as this is not entirely unknown in the literature,
see Rachev~\cite[eq.~(3.2.22)]{Rach91}.

There is an alternative definition, 
due to Strassen~\cite[Cor.~to Thm~11]{Str65},
of the Prohorov metric in terms of an optimal coupling.
The (parameterised) {\it Ky Fan distance\/} 
$K_\slack(X,Y)$ between random variables (r.v's) 
$X$, $Y$ on $\states$ is defined as
$$
K_\slack(X,Y)\defeq\inf\big\{\eps:\Pr[d(X,Y)>\slack\eps]\leq\eps\big\}.
$$
Denote by $\law(X)$ the law (distribution) of r.v.~$X$.
\begin{theorem}\label{thm:Strassen}
Suppose $\pi$ and $\pi'$ are probability distributions on~$\states$.
Then $\Prohp(\pi,\pi')$ is the infimum of $K_\slack(X,Y)$ over all
pairs $(X,Y)$ of coupled $\states$-valued r.v's such that 
$\law(X)=\pi$ and $\law(Y)=\pi'$.
\end{theorem}
(The theorem in this form is from Garc\'\i a-Palomares 
and Gin\'e~\cite{G-PG77}.)

\begin{remark}
Strassen states Theorem~\ref{thm:Strassen} for the case $\slack=1$,
but the proof clearly holds for arbitrary $\slack\geq0$.
(To avoid delving into the proof, one could simply scale
the metric~$d$.)
The case $\slack=0$ is the well-known Optimal Coupling Theorem.
It is in the proof of Theorem~\ref{thm:Strassen} that the technical 
assumption of separability is used.
\end{remark}

%

One last definition, and we'll be ready to formalise
conditions~(1)--(3).  The {\it total variation 
distance\/} between two measures $\pi$ and $\pi'$ 
on~$\states$ is 
$$
\|\pi-\pi'\|_\mathrm{TV}\defeq\Proh_0(\pi,\pi')=
   \inf\big\{|\pi(A)-\pi'(A)|:A\text{ closed, }A\subseteq\states\big\}.
$$
The the {\it variation threshold time\/}~\cite[\S4.3]{AldFill}
of the Markov chain $(\states,P)$ is defined to be
$$
\tau_1\defeq\min\big\{t: \|P^t(x,\cdot)-P^t(x',\cdot)\|_\mathrm{TV}\leq e^{-1},
\text{ for all $x,x'\in\states$}\big\}.
$$
The choice of threshold $e^{-1}$ is somewhat arbitrary.
There are other, slightly different notions of $\ell_1$ mixing,
but they are equivalent for our purposes.  
In algorithmic applications, one often estimates the 
probability $\pi(A)$ of some event~$A$ in the stationary 
distribution by taking a suitably sized sample from 
the $t$-step distribution $P^t(x,\cdot)$.  There are two
sources of error in this process:  the sampling error, 
and the error occasioned by using 
$P^t(x,\cdot)$ in place of~$\pi(\cdot)$.  
The  variation threshold time is important precisely because it 
is a worst-case bound on the latter. 

\section{Main result}
Now all the definition are in place we can state the main result.
\begin{theorem}\label{thm:main}
Suppose for some $\slack,C,\delta\geq0$:  
\begin{enumerate}
\item $\Prohp(\PP(x,\cdot),P(x,\cdot))\leq \delta$, 
for all $x\in\statesP$;  
\item $\Prohp(P(x,\cdot),P(x',\cdot))\leq C\,d(x,x')$, for all 
$x,x'\in\states$;
\item 
The Markov chain defined by $P$ is ergodic, 
with stationary distribution~$\pi$,
and variation threshold time~$\tau_1$.
\end{enumerate}
Then $\Prohp(\PP^t(x,\cdot),\pi)\leq\eps$ provided
$t\geq t_\eps\defeq \lceil\ln(2e/\eps)\tau_1\rceil$ and, additionally:
\begin{itemize}
\item in the case $\slack C< 1$, 
$$\delta\leq\frac{(1-\slack C)\eps}{2t_\eps};$$
\item in the case $\slack C=1$, 
$$\delta\leq\frac{\eps}{t_\eps(t_\eps+1)};$$ 
\item and in the case $\slack C>1$,
$$\delta\leq\frac{(\slack C-1)^2\eps}{2(\slack C)^{t_\eps+1}}.$$
\end{itemize}
\end{theorem}

\begin{remarks}
\begin{itemize}
\item The key point is that if $\slack C<1$ 
then $\PP$ does not need to approximate~$P$
to excessive accuracy, but only to within $O(\tau_1^{-1})$.
In contrast, when $\slack C>1$, the required accuracy scales 
exponentially with~$\tau_1$.  So, for example, real arithmetic would have 
to be carried out to a number of significant digits scaling 
linearly with~$\tau_1$.  
In the boundary situation, $\slack C=1$, the  
required accuracy scales as $O(\tau_1^{-2})$.

\item All three behaviours described in Theorem~\ref{thm:main}
actually occur, and are not artifacts of the proof.  Examples 
will be provided in \S\ref{sec:ex}.

\item In \S\ref{sec:app}
we shall see that the ball-walk, at least of the lazy kind,
fits the most favourable case, $\slack C<1$.

\item We can recover something akin to one of 
Azar et al.'s results~\cite{ABKLP96} by
setting $\slack=0$, $C=1$ and $d$ to be the discrete metric.
Observe that condition~2 of the theorem becomes vacuous,
and $\Proh_0$ is just total variation distance.
Note that Azar et al.\null{} express their condition~3
in terms of $\ell_2$~mixing time (spectral gap).
\end{itemize}
\end{remarks}

\begin{proof}[Proof of Theorem~\ref{thm:main}]
Set $t=t_\eps=\lceil\ln(2e/\eps)\tau_1\rceil$.  
Let $(X_i)$ and $(\Xhat_i)$ be Markov chains with
transition kernels $P$ and~$\PP$, respectively, 
starting at a fixed state $X_0=\Xhat_0=a\in\statesP$.
Note that $t$ has been chosen so that
$\Prohp(\law(X_t),\pi)\leq\|\law(X_t)-\pi\|_{\mathrm{TV}}\leq\eps/2$.
(See, e.g., Aldous and Fill~\cite[\S4, Lemma~5]{AldFill}.)

We'll couple $(\Xhat_i)$ and $(X_i)$ so that 
\begin{equation}\label{eq:XhatvsX}
\Prohp(\law(\Xhat_t),\law(X_t))\leq\eps/2.
\end{equation}
This will be possible provided
$\delta$ satisfies the appropriate condition laid down
in the statement of Theorem~\ref{thm:main}.  
To see this, let $D_i\defeq d(\Xhat_i,X_i)$ denote the 
divergence of the two Markov chains at time~$i$. 
Consider the situation at time $i-1$.
Suppose we have constructed a realisation
of the coupled process 
$$(a,a)=(\Xhat_0,X_0), (\Xhat_1,X_1)\ldots 
(\Xhat_{i-1},X_{i-1})=(\bhat,b).$$
Conditioned on $(\Xhat_{i-1},X_{i-1})=(\bhat,b)$
we have 
\begin{align*}
\Prohp(\law(\Xhat_i),\law(X_i))&\leq 
\Prohp(\law(\Xhat_i),P(\bhat,\cdot))+\Prohp(P(\bhat,\cdot),\law(X_i))\\
&=\Prohp(\PP(\bhat,\cdot),P(\bhat,\cdot))
   +\Prohp(P(\bhat,\cdot),P(b,\cdot))\\
&\leq\delta+CD_{i-1},
\end{align*}
where the final inequality uses conditions (1) and~(2) of the theorem.
According to Theorem~\ref{thm:Strassen}, 
we may couple $\Xhat_i$ and~$X_i$ so that
$$
\Pr\big[D_i>\slack(CD_{i-1}+\delta)\big]\leq CD_{i-1}+\delta.
$$
Iterating this construction, it follows, by induction on~$i$,
that   
\begin{equation}\label{eq:recSoln}
\Pr\left[D_t>\slack\delta\sum_{i=0}^{t-1}(\slack C)^i\right]\leq 
   \delta\sum_{i=0}^{t-1}(t-i)(\slack C)^i;
\end{equation}

Considering first the case $\slack C<1$, 
we may sum the series in~(\ref{eq:recSoln}) 
to obtain
\begin{equation}\label{eq:recSoln'}
\Pr\left[D_t>\frac{\slack\delta(1-(\slack C)^t)}{1-\slack C}\right]\leq 
   \frac{\delta t}{1-\slack C}
   -\frac{\slack\delta C(1-(\slack C)^t)}{(1-\slack C)^2},
\end{equation}
which entails
\begin{equation}
\Pr\left[D_t>\frac{\slack\delta}{1-\slack C}\right]\leq 
   \frac{\delta t}{1-\slack C}.
\end{equation}
Our goal is to attain
\begin{equation}\label{eq:goal}
\Pr[D_t>\slack\eps/2]\leq\eps/2,
\end{equation}
since this 
implies inequality~(\ref{eq:XhatvsX}) 
through Theorem~\ref{thm:Strassen}.
The analysis of the case $\slack C<1$
is completed by noting that to achieve the goal 
it is sufficient that $\delta\leq(1-\slack C)\eps/2t$.

Now turn to the case $\slack C=1$.  Summing the series 
in~(\ref{eq:recSoln}) in this case yields
$$
\Pr[D_t>\slack\delta t]\leq\frac{\delta t(t+1)}2.
$$  
We achieve (\ref{eq:goal}) provided
$\delta\leq\eps/t(t+1)$.

The final case, $\slack C>1$ is handled in a very similar 
manner to the first.  In this case we find
\begin{equation}
\Pr\left[D_t>\frac{\slack\delta((\slack C)^t-1)}{\slack C-1}\right]\leq 
   \frac{\slack\delta C((\slack C)^t-1)}{(\slack C-1)^2},
\end{equation}
and that (\ref{eq:XhatvsX}) is achieved provided
$$
\delta\leq\frac{(\slack C-1)^2\eps}{2(\slack C)^{t+1}}.
$$

In conclusion, we have shown that
$$
\Prohp(\law(\Xhat_t),\pi)\leq  \Prohp(\law(\Xhat_t),\law(X_t))+
\Prohp(\law(X_t),\pi) \leq \eps/2+\eps/2=\eps,
$$
as required.
If $t>t_\eps$, we simply delay starting the coupling until $t_\eps$ 
steps from the end.
\end{proof}

\section{Counterexamples}\label{sec:ex}
We demonstrate in this section
that the dependence on~$\tau_1$ indicated by Theorem~\ref{thm:main}
is correct:  i.e., linear in the case $\slack C<1$,
exponential in the case $\slack C>1$, and quadratic at the 
boundary.

In applications we are thinking mainly of uncountable state spaces.
However, for convenience, the counterexamples 
will all be finite Markov chains.

\subsection{``Convergent'' case}
The heading is intended to indicate the case $\slack C<0$.
We'll set $\slack=0$ (i.e., our measure of convergence 
is total variation distance) and $C=1$, though the construction
would work equally well for a range of $\slack, C$ satisfying 
$\slack C<1$. 

The state space in this counterexample
is $\states\defeq\{\omega_j:0\leq j<n\}$.
Identify the state $\omega_j$ with the point 
$(n\cos(2j\pi/n), n\sin(2j\pi/n))$ in $\reals^2$,
so that the states are equally spaced points around 
a circle of radius~$n$.
The metric~$d$ is just Euclidean distance.

Define transition probabilities for the Markov chain 
(from state~$\omega_j$) according to the following trial:
\begin{itemize}
\item With probability $1/n$, set $j':=0$.
\item Otherwise (with probability $1-1/n$), set $j':=(j+1)\bmod n$.
\end{itemize}
The new state is $\omega_{j'}$.  Informally, we 
move relentlessly clockwise around the circle, except that with 
probability $1/n$ we perform a ``reset'' and return to distinguished
vertex $\omega_0$. 
Since $C=1$ and the Euclidean distance 
between any pair of states is at least~$1$,
condition~(2) of Theorem~\ref{thm:main} is vacuously true.

It is easy to verify, by coupling, 
that the variation threshold 
time~$\tau_1$ is $O(n)$.  
Simply take two copies of the Markov chain 
and couple the resets.  A synchronised reset occurs 
within $n$~steps with probability at least $1-e^{-1}$,
so $\tau_1\leq n$.  (See, e.g., Aldous~\cite[Lemma~3.6]{Ald83}.)

Define $\PP$ as $P$ but with reset probability~$4/n$ in place 
of~$1/n$.  We claim that with 
$\PP$ there is significantly 
lower probability of observing $j\geq n/2$.
Thus the stationary distributions are quite far apart 
in total variation distance (which is Prohorov metric
with parameter $\slack=0$).

The justification of this claim runs as follows.  Assume
for convenience that $n$ is even,
and fix a time step $t\geq n$.
The probability that we observe $j\geq n/2$ (in the 
$P$~version) is at least
$$
\Pr(\text{no reset in past $n/2$ steps}\wedge 
\text{at least one reset in past $n$ steps})\\
$$
This for large $n$ is close to 
$e^{-1/2}(1-e^{-1/2})\geq 0{\cdot}238$.
In contrast, for the $\PP$ version,
the probability that we observe $j\geq n/2$ is at most 
$$
\Pr(\text{no reset in past $n/2$ steps}),
$$
which for large $n$ is close to 
$e^{-2}\leq 0{\cdot}136$.  Comparing with previous bound,
it will be seen that the 
two stationary distributions differ by at least $0{\cdot}1$
in total variation distance.

So we certainly need to insist on $\delta<4/n-1/n=3/n$
if we want to guarantee that the 
stationary distributions of the two Markov chains 
are closer than $\eps=0{\cdot}1$ in variation distance.  
In particular, we could not replace the $\tau_1$~factor in the 
first case of Theorem~\ref{thm:main} by anything growing more slowly.

\subsection{``Neutral'' case}
This is the boundary case $\slack C=1$.
The state space for this example is
$\states\defeq\{\omega_{i,j}:\allowbreak 0\leq i,j<n\}$.
Identify state $\omega_{i,j}$ with the point 
$(n\cos(2j\pi/n),\allowbreak n\sin(2j\pi/n),\allowbreak 5i/n^2)$ 
in $\reals^3$.
(There are $n$ circles, in $n$~layers closely packed 
the $z$-dimension, each containing $n$ evenly spaced states.)
The metric~$d$ is again Euclidean distance.
Define
$$
r(i)\defeq \begin{cases}1/n&\text{if $i<n/5$;}\\
   5i/n^2&\text{if $n/5\leq i<4n/5$;}\\
   4/n&\text{if $i\geq 4n/5$.}\end{cases}
$$

Define transitions probabilities (from state~$\omega_{i,j}$)
according to the following trial:
\begin{enumerate}
\item 
\begin{itemize}
\item With probability $r(i)$, set $j':=0$;
\item Otherwise (with probability $1-r(i)$), set $j':= j+1\mod n$.
\end{itemize}

\item 
\begin{itemize}
\item With probability $2/3$, set $i':=\max\{i-1,0\}$;
\item Otherwise (with probability $1/3$), set $i':=\min\{i+1,n-1\}$.
\end{itemize}
\end{enumerate}
The new state is $\omega_{i',j'}$.
Informally:  owing to the drift in the $z$-dimension, 
we quickly gravitate to $i=0$ layer and stay close to it.  
Within the layer,
we move clockwise around the cycle, except that with 
probability $r(i)$ we perform a 
``reset'' and return to one of the 
distinguished states $\omega_{i,0}$.

We set $\slack=1$, in other words we measure convergence
in the standard Prohorov metric.
It is routine to verify that 
$\Proh\big(P(x,\cdot),P(x',\cdot)\big)\leq d(x,x')$,
so that $C=1$, and we are in the $\slack C=1$ (boundary) regime.  
(We need only check pairs of states of the form
$(x,x')=(\omega_{i,j},\omega_{i',j})$, 
i.e., pairs which agree in their second index,
since other pairs of states have separation $d(x,x')>1$.
Indeed, by the triangle inequality, we need only
check pairs of the form 
$(x,x')=(\omega_{i,j},\omega_{i+1,j})$.
There is a natural coupling of transitions 
from these adjacent states $x$ and~$x'$ 
such that the new states are 
within distance $5/n^2$ of each other with probability at 
least $1-5/n^2$.)

As before, we can show that $\tau_1=O(n)$
using a coupling argument.  Consider two copies of the 
Markov chain started in different states.
In the first phase, couple on~$i$ using the identity coupling.
Coupling (of the $i$-index) 
occurs at or before the first occasion at which both 
copies have visited $i=0$ layer.  After this point the coupled 
versions always agree as to the level.  This happens 
with high probability within $4n$ steps.
In the second phase, we couple on~$j$.
We do the natural thing 
and synchronise the resets (just as in the case $\slack C<1$).  
Again, we can arrange for a synchronised reset within
$2n$ steps with high probability.

Define $\PP$ as $P$ but with drift on $i$ reversed.
The intuition is that we quickly gravitate to layer $i=n-1$ and 
remain close to it.  We then circle as before, but with 
much higher reset probability.  We claim, as before, that with 
$\PP$ there is a lower probability of observing $j\geq n/2$.
Thus the stationary distributions are quite far apart 
in Prohorov metric.

The justification of this claim runs as follows.
Denote by $\calE$ the event (in the $P$~Markov chain)
\begin{quote}
``in the previous $n$ steps, $i$ has remained in range $[0,n/5]$''
\end{quote}
Fix a time step $t\geq 6n$.
The probability that we observe $j\geq n/2$ (in the $P$~version) 
is at least
{\setlength{\multlinegap}{0pt}
\begin{multline*}
\Pr(\calE\wedge\text{no reset in past $n/2$ steps}\wedge 
\text{at least one reset in past $n$ steps})\\
 \null=\Pr(\calE)\Pr(\text{no reset in past $n/2$ steps}
 \wedge\text{at least one reset in past $n$ steps}\mid\calE).
\end{multline*}}
This for large $n$ is close to 
$e^{-1/2}(1-e^{-1/2})\geq 0{\cdot}238$.
In contrast, for the $\PP$ version, denote by $\calE'$ the event
\begin{quote}
``in the previous $n$ steps, $i$ has remained in range $[4n/5,n]$'',
\end{quote}
the probability that we observe $j\geq n/2$ is at most 
\begin{multline*}
\Pr(\neg\calE'\vee\text{no reset in past $n/2$ steps})\\
 \null\leq\Pr(\neg\calE')+\Pr(\text{no reset in past $n/2$ steps}\mid \calE').
\end{multline*}
The latter probability for large~$n$ is close to 
$e^{-2}\leq 0{\cdot}136$.  Comparing with the previous
estimate, we see the 
two stationary distributions differ by at least $0{\cdot}1$
in the Prohorov metric.

So we certainly need to insist on $\delta <10n^{-2}$
to bring the stationary distributions of the two Markov chains
within $\eps=0{\cdot}1$ in the Prohorov metric.
In particular, we could not, e.g., replace exponent~2
in the second case of Theorem~\ref{thm:main} by anything smaller.

\subsection{``Divergent'' case}
The state space here is $\states\defeq\{2^{-n}i:0\leq i<2^n\}$, 
and the metric $d:\states^2\to\reals^+$ is given by $d(x,y)\defeq|x-y|$.
Define the function $G:\states\to\states$ by
$$
G(x)\defeq\begin{cases}
2x,&\text{if $x<1/2$;}\\
2(1-2^{-n}-x)&\text{if $x\geq1/2$.}
\end{cases}
$$
(This is the Lorenz ``tent map'' 
of dynamical systems~\cite[eq.~(2.5.2)]{SH96}, 
adapted to the discrete situation.)  Then the 
transition kernel 
$$
P(x,y)\defeq
\begin{cases}1/2,&\text{if $y\in\{G(x),G(x)+2^{-n}\}$;}\\
0,&\text{otherwise,}
\end{cases}
$$
defines an ergodic Markov chain with 
stationary distribution~$\pi$ uniform on~$\states$. 
Why is this?  View 
$x=0{\cdot}x_{1}x_{2}\ldots x_{n}\in\states$ as
an $n$-bit binary fraction.  
Then 
$$
G(x)=\begin{cases}
   0{\cdot}x_{2}x_{3}\ldots x_{n}0,&\text{if $x<1/2$;}\\
   0{\cdot}\xbar_{2}\xbar_{3}\ldots \xbar_{n}0,&\text{if $x\geq1/2$,}
\end{cases}
$$
where $\xbar_i\defeq1-x_i$.  That is to say, $G$ can be viewed
as a left shift, followed (possibly) by complementation.
(C.f.\ two's complement arithmetic.)
So one step of the Markov chain can be viewed as a left shift,
followed (possibly) by complementation, and concluded
by appending a random bit.
Thus $P^t(x)$ for any $x\in\states$ and $t\geq n$ is a 
binary fraction formed of independent, symmetric
Bernoulli r.v's.  
We see from this argument that $\tau_1=n$. (Notice that 
distance from stationarity drops from $1/2$ to~$0$ 
between time $t=n-1$ and time $t=n$!)
Set $\slack=1$, and 
observe that $\Proh(P(x,\cdot),P(x',\cdot))\leq 2d(x,x')$,
so that we are $\slack C>1$ regime.
(In light of the triangle inequality, we just need to check
pairs $(x,x')$ with $x'=x+2^{-n}$.)

Now define an approximating Markov chain:
$$
\PP(x,y)\defeq
\begin{cases}
   3/4,&\text{if $y=G(x)$ and $2^{n-1}G(x)$ is even;}\\
   1/4,&\text{if $y=G(x)$ and $2^{n-1}G(x)$ is odd;}\\
   3/4,&\text{if $y=G(x)+2^{-n}$ and $2^{n-1}G(x)$ is odd;}\\
   1/4,&\text{if $y=G(x)+2^{-n}$ and $2^{n-1}G(x)$ is even;}\\
   0,&\text{otherwise.}
\end{cases}
$$
Note that $\Proh(\PP(x,\cdot),P(x,\cdot))\leq2^{-n}$.

The interpretation of the Markov chain 
defined by $\PP$ in terms of binary 
fractions is similar to before, only now the random bit 
appended is with probability $3/4$ equal to the bit 
immediately to its left.  So, for any $t\geq n$,
$\PP^t(x,A)=3/4$, where $A=\states\cap\big([0,1/4)\cup[3/4,1)\big)$.
In contrast, 
$P^t(x,A')=2/3$, where $A'=\states\cap\big([0,1/3)\cup[2/3,1)\big)$.
Now $A'\supseteq A^\eps$ with $\eps=1/12$.
Thus 
$$
\Proh(\PP^t(x,\cdot),\pi)=\Proh(\PP^t(x,\cdot),P^t(x,\cdot))\geq1/12,
$$
where $\pi$ is the stationary (uniform) distribution.
The bottom line is that the transition kernels 
$\PP$ and~$P$ are very close and 
variation threshold time is short, 
but  that the stationary distributions 
of the two Markov chains are nevertheless
far apart.
The exponential dependence of $\delta$ on $\tau_1$ 
in the third case of Theorem~\ref{thm:main}
is unavoidable.

\section{Application:  ball walk of Lov\'asz and Simonovits}\label{sec:app}
Recall the ball walk of Lov\'asz and Simonovits~\cite{LS93}
in its ``lazy'' version.  The situation is as follows.
$K\subset\reals^n$ is a convex body in $n$-dimensional Euclidean space.  
For $x\in\reals^n$ and $r\in\reals^+$, 
$B_n(x,r)$ denotes the $n$-dimensional (closed) ball centred 
at~$x$.  Procedurally, the lazy walk $(X_t:t\in\nats)$ is described 
by the following trial (where the current state is $X_t=x\in\reals^+$):
\begin{enumerate}
\item Choose $y\in B_n(x,r)$ u.a.r.
\item If $y\in K$ then $X_{t+1}:= y$ else $X_{t+1}:= x$.
\end{enumerate}
Alternatively, the transition kernel is 
\begin{equation}\label{eq:BWkernel}
P(x,A)\defeq 
\begin{cases}
   \mu_n\big(B_n(x,r)\cap (\Kbar \cup A)\big)/v_n(r),&\text{if $x\in A$;}\\
   \mu_n(B_n(x,r)\cap K\cap A)/v_n(r),&\text{otherwise},
\end{cases}
\end{equation}
where $\mu_n$ is Lebesgue measure, $\Kbar$ denotes the complement of~$K$,
and $v_n(r)\defeq\mu_n(B_n(0,r))$ the 
volume of the $n$-dimensional ball of radius~$r$.

To apply Theorem~\ref{thm:main}, we want to find a constant~$C$
such that 
$$
\Prohp\big(P(x,\cdot),P(x',\cdot)\big)\leq C\,d(x,x')= C\,\|x-x'\|_2,
$$
since $d$ is here Euclidean distance.
For this part of the calculation the value of~$\slack$
is immaterial (even $\slack=0$ will do), so we'll 
defer the choice of~$\slack$ until later. 

We could work directly from~(\ref{eq:BWkernel}),
but it seems easier to go via Theorem~\ref{thm:Strassen}.
Let $d=\|x-x'\|_2$.  For convenience, 
let $x=d\unit/2$ and $x'=-d\unit/2$,
where $\unit$ is the unit vector parallel 
to the first coordinate axis.
Define a coupling $(Y,Y')$ with  
$\law(Y)=P(x,\cdot)$ and $\law(Y')=P(x',\cdot)$
according to the trial
\begin{enumerate}
\item Choose $y\in B_n(x,r)$ u.a.r.
\item If $y\in B_n(x',r)$ then $y':= y$ else $y':=\ybar$,
where $\ybar$ is the reflection of~$y$ in the plane $\xi\cdot\unit=0$.
\item 
\begin{itemize}
\item If $y\in K$ then $Y:= y$ else $Y:= x$;
\item If $y'\in K$ then $Y':= y'$ else $Y':= x'$.
\end{itemize}
\end{enumerate}
Note that $Y=Y'$ unless $y\in B_n(x,r)\setminus B_n(x',r)$.
Now $\mu_n(B_n(x,r)\setminus B_n(x',r))$ is bounded above
by the volume of a $n$-dimensional cylinder with height~$d$
and cross-sectional $(n-1)$-dimensional volume $v_{n-1}(r)$.  
Thus $\Pr(Y\not=Y')\leq d\,v_{n-1}(r)/v_n(r)$, and hence 
$$
\Prohp(P(x,\cdot),P(x',\cdot))\leq \frac{v_{n-1}(r)}{v_n(r)}\,d\leq Cd
$$
where $C=\Theta(\sqrt n/r)$.  (Note that the inequality 
holds for any $\slack$, even $\slack=0$.)
By setting $\slack=1/2C=\Theta(r/\sqrt n\,)$ 
we place ourselves in the first
(most favourable) case of Theorem~\ref{thm:main}.\footnote{In 
applications of the ball walk, the radius $r$ is typically 
of order $1/\sqrt n$, so that $\slack=\Theta(1/n)$.}

Now, under the simplifying assumption that
the convex body~$K$ does not have sharp corners,
the variation threshold time is
\begin{equation}\label{eq:vttBW}
\tau_1=O\left(\frac{D^2n^2\ln(D/r)}{r^2}\right),
\end{equation}
where $D$ is the diameter of~$K$.

\begin{remark}
See \cite[Thm~6.7 and Cor.~6.8]{Jer03} for 
more detail, including a precise explanation of the 
requirement of having no ``sharp corners''.
Note that the radius of the ball defining the ball 
walk is usually denoted~$\delta$;
we have used~$r$ instead to avoid a notational clash.  
For general convex bodies~$K$, the mixing time is essentially
as given in~(\ref{eq:vttBW}), but one has to take care 
over the distribution of the start state of the walk,
since the ball walk in its lazy variant may get trapped
for long periods near points on the boundary of~$K$ 
of tight curvature.
\end{remark}

From the above considerations, it can be seen that
the transition kernel~$\PP$ 
of the ball walk as implemented is not required to 
approximate the ideal transition kernel very closely;
specifically we require, according to Theorem~2, 
$\Prohp(\PP(x,\cdot),P(x,\cdot))\leq \delta$,
where 
\begin{equation}\label{eq:deltBW}
\delta = O\left(\frac{\eps r^2}{D^2n^2\ln(D/r)}\right).
\end{equation}
This is consistent 
with Lov\'asz and Simonovits's observation that 
for their algorithm real numbers need only be carried to 
$O(\log n)$ digits.

Some concise notes on how to achieve~(\ref{eq:deltBW}).
Assume, as a starting point, procedures that sample points 
from distributions that are close to $N(0,1)$
(Gaussian with mean~0 and variance~1)
and to $U(0,1)$ (uniform on $[0,1]$).
A standard approach to sampling a point u.a.r.\null{} 
from $B_n(0,r)$ is the following (with step~2 omitted):
\begin{enumerate}
\item Let $\Nrv_1,\Nrv_2,\ldots,\Nrv_n$ be i.i.d.\null{} samples
from $N(0,1)$.
\item If $R=\sqrt{\Nrv_1^2+\Nrv_2^2+\cdots+\Nrv_n^2}<\frac12\sqrt n$ 
declare the trial void and start again at step~1.
\item Set $S=r\,(\Nrv_1,\Nrv_2,\ldots,\Nrv_n)/R$.
\item Let $U$ be a sample from the uniform distribution on $[0,1]$,
and return $W=U^{1/n}S$. 
\end{enumerate}
We assume throughout that arithmetic is exact, in order 
to focus on sampling errors.
Assume that $\Nrv_1,\Nrv_2,\ldots,\Nrv_n$ and $U$ are sampled 
perfectly from distributions $N(0,1)$ and $U(0,1)$.
Then $(\Nrv_1,\Nrv_2,\ldots,\Nrv_n)$ is distributed according to 
an $n$-dimensional symmetric Gaussian distribution, 
and is in particular rotationally symmetric.
Thus, with or without step~2, $S$ is distributed uniformly
over the surface of $B_n(0,r)$.  The finally step spreads 
the distribution uniformly into the interior of $B_n(0,r)$.
The unusual step~2 is included to avoid a small error being 
blown up in the unlikely event that~$R$ is close to~0.

Without loss of generality, assume that ball walk is
at the origin at time step~0. Its location~$Y$ at time step~1  
is obtained by applying the rejection rule to the r.v.~$W$;
explicitly, $Y=W$ if $W\in K$, and $Y=0$ otherwise.
Now suppose that we have only approximations $\Nrvhat_i$
and~$\Uhat$ to the perfect samples.  Specifically, 
suppose 
\begin{equation}\label{eq:NandU}
\Proh\big(\law(\Nrvhat_i),N(0,1)\big)=O(\delta/n)\quad\hbox{and}\quad 
\Proh\big(\law(\Uhat),U(0,1)\big)=O(\delta^2),
\end{equation} 
where $\delta$, given by~(\ref{eq:deltBW}), 
is the deviation we are prepared to tolerate in $\law(\Yhat)$,
the approximate version of~$Y$.
(Specifically, we are aiming at
$\Prohp(\law(\Yhat),\law(Y))\leq\delta$.)
Suppose that we run through the above trial, replacing 
the perfectly distributed r.v's by their hatted,
imperfect approximations $\Shat$, $\What$ and finally~$\Yhat$,
which arises from the rejection rule:
$\Yhat=\What$ if $\What\in K$, and $\Yhat=0$ otherwise.

Now couple the hatted and 
unhatted r.v's as suggested by Theorem~\ref{thm:Strassen}.
The build-up of errors is summarised in the following 
table.  The penultimate row relates to the approximate 
proposal move~$\What$ for the ball walk, sampled
according to the four-step trial described earlier,
and the final row to the result of applying the rejection rule.
The interpretation of (say) the third line of the table 
is that we may couple $S$ and~$\Shat$ so that  
$\|S-\Shat\|_2=O(r\delta/n)$ with probability $1-O(\delta)$.

{\renewcommand{\arraystretch}{1.4}
$$
\begin{tabular}[h]{|c|c|c|}
\hline
Random variable~$\Xhat$& $\|X-\Xhat\|_2$ bounded by$\ldots$& except with probability$\ldots$\\
\hline\hline
$\Nrvhat_i$ & $O(\delta/n)$ & $O(\delta/n)$\\ 
\hline
$(\Nrvhat_1,\Nrvhat_2,\ldots,\Nrvhat_n)$ & $O(\delta/\sqrt n\,)$ & $O(\delta)$\\
\hline
$\Shat$& $O(r\delta/n)$ & $O(\delta)+O(\delta/n)= O(\delta)$\\
\hline
$\What$& $O(r\delta/n)$ & $O(\delta)+O(\delta)= O(\delta)$\\
\hline
$\Yhat$& $O(r\delta/n)$ & $O(\delta)+O(\delta)= O(\delta)$\\
\hline
\end{tabular}
$$}

The rows of the table may be checked as follows.
Row~2 is straightforward.  In row~3 we need to 
be concerned about the trial being declared void in the hatted 
trial and not in the unhatted, or vice versa.  For this to
occur, $R$ must be within $O(\delta/\sqrt n\,)$ of $\frac12\sqrt n$,
an event whose probability may be (crudely) 
bounded by $O(\delta/\sqrt n\,)\times O(1/\sqrt n\,)=O(\delta/n)$.
(The density of the r.v.~$R$ is unimodal, and 
achieves its maximum at the point $\sqrt{n-1}$;
so the density of~$R$ at $\frac12\sqrt n$ can be at 
most $O(1/\sqrt n\,)$.)
In row~4, we need to be concerned about errors being 
magnified when $U$ is close to~0.  We deal with this simply 
by giving everything away if $\Uhat=O(\delta)$.
In the final row, our concern is with the 
event $\Yhat\in K$ and $Y\not\in K$ (or vice versa).
For this event, we must have $Y\in (K^\eta\setminus K)\cap B(0,r)$,
where, as usual, $K^\eta$ denotes the Minkowski sum of~$K$ and a ball 
of radius~$\eta$, and $\eta=O(r\delta/n)$.  
Now 
$$
\mu_n\big((K^\eta\setminus K)\cap B(0,r)\big)\leq 
\mu_n\big(B(0,r)^\eta \setminus B(0,r)\big),
$$
and so 
$$
\frac{\mu_n\big((K^\eta\setminus K)\cap B(0,r)\big)}{\mu_n(B(0,r))}\leq
\frac{\mu_n\big(B(0,r)^\eta \setminus B(0,r)\big)}{\mu_n(B(0,r))}=O(\delta).
$$

Recall that we have set $\lambda=\Theta(r/\sqrt n\,)$,
from which it follows that $r\delta/n=O(\slack\delta)$.
In summary, then, to obtain a close approximation
to the ball-walk it is enough that the various
samples from the Gaussian and uniform distributions
satisfy~(\ref{eq:NandU}), 
where $\delta$ is given by~(\ref{eq:deltBW}). 

\begin{remark}
It is unlikely that one would want, 
in the analysis of a new algorithm,
to repeat a calculation such as the one 
given above in a similar level of detail.  Nevertheless, 
it would be comforting to verify, in practical 
situations, that one was working in one 
of the two favourable cases in Theorem~\ref{thm:main}:  
it would then follow by
more informal reasoning that logarithmic 
(number of bits or significant digits)
accuracy would suffice. 
\end{remark}
\bibliography{draft}

\end{document}